\documentclass[reqno,14pt,centertags]{amsart}
\usepackage{amsmath,amsthm,amscd,amssymb,latexsym,upref,stmaryrd}

\usepackage{hyperref}
\newcommand*{\mailto}[1]{\href{mailto:#1}{\nolinkurl{#1}}}

\newcommand{\B}{$\hfill\Box$}
\newcommand{\al}{\alpha}

\newcommand{\de}{\delta}
\newcommand{\la}{\lambda}
\newcommand{\om}{\omega}

\newcommand{\vv}{\varphi}
\newcommand{\iy}{\infty}

\theoremstyle{definition}

\begin{document}
\title[Recovering Nonlocal Differential Pencils]
{Recovering Nonlocal Differential Pencils}

\author[C.~F.~Yang]{Chuan-Fu Yang}
\address{Department of Applied Mathematics, Nanjing University of
Science and Technology, Nanjing, 210094, Jiangsu, People's
Republic of China}
\email{\mailto{chuanfuyang@njust.edu.cn}}

\author[V.~Yurko]{Vjacheslav Yurko}
\address{Department of Mathematics, Saratov University,
Astrakhanskaya 83, Saratov 410012, Russia}
\email{\mailto{yurkova@info.sgu.ru}}

\date{\today}

\begin{abstract}
Inverse problems for differential pencils with nonlocal conditions are investigated.
Several uniqueness theorems of inverse problems from the Weyl-type function and spectra
are proved, which are generalizations of the well-known Weyl function and Borg's inverse
problem for the classical Sturm-Liouville operator.
\end{abstract}

\maketitle

\section{Introduction}

Problems with nonlocal conditions arise in various fields
of mathematical physics \cite{d,g,i,1,yin}, biology and biotechnology \cite{n,s}, and
in other fields. Nonlocal conditions come up when
value of the function on the boundary is connected to values inside the domain.
Recently problems with nonlocal conditions are paid much attention for them
in the literature.

In this paper we study inverse spectral problems for differential pencils
\begin{align}\label{1}
y''(x)+[\lambda^2-2\lambda p(x)-q(x)]y(x)=0,\quad x\in (0,T),
\end{align}
and with nonlocal linear conditions
\begin{align}\label{2}
U_j(y):=\int_0^T y(t)d\sigma_j(t)=0,\quad j=1,2.
\end{align}
Here $p\in AC[0,T]$ (absolutely continuous function) and $q\in L(0,T)$ are complex-valued functions, $\sigma_j(t)$ are
complex-valued functions of bounded variations and are continuous from
the right for $t>0.$ There exist finite limits $H_j:=\sigma_j(+0)-
\sigma_j(0).$ Linear forms $U_j(y)$ in (\ref{2}) can be written as forms
\begin{align}\label{3}
U_j(y):=H_j y(0)+\int_0^T y(t)d\sigma_{j0}(t),\quad j=1,2,
\end{align}
where $\sigma_{j0}(t)$ in (\ref{3}) are complex-valued functions of bounded
variations and are continuous from the right for $t\ge 0$, and $H_1\ne 0.$

A complex number $\lambda_0$ is called an eigenvalue of the problem (\ref{1}) and (\ref{2}) if
equation (\ref{1}) with $\lambda=\lambda_0$ has a nontrivial solution $y_0(x)$
satisfying conditions (\ref{2}); then $y_0(x)$ is called the eigenfunction of the problem
(\ref{1}) and (\ref{2}) corresponding to the eigenvalue $\lambda_0$. The number of linearly independent solutions of the
problem (\ref{1}) and (\ref{2}) for a given eigenvalue $\lambda_0$ is called the multiplicity of $\lambda_0$.

Classical inverse problems for
Eq.(\ref{1}) with two-point boundary conditions have been studied
fairly completely in many works (see the monographs \cite{2,3,4,5,Yurko2,GGU} and the
references therein). The theory of nonlocal inverse spectral problems
now is only at the beginning because of its complexity. Results
of the inverse problem for various nonlocal operators can be
found in \cite{6,7,8,9,10,11,12,13,14,15}.

In this work by using Yurko's ideas of the method of spectral mappings \cite{5} we prove uniqueness theorems for the
solution of the inverse spectral problems for Eq.(\ref{1}) with nonlocal
conditions (\ref{2}). In Section 2 we formulate our main results (Theorems 1 and 2). Section 3
introduces some properties of spectral characteristics.
The proofs of Theorems 1 and 2 are given in Section 4. In Section 5 we
provide two counterexamples related to the statements of the inverse problems
(see also \cite{8,yy}). 
In Section 6, as a consequence of Theorem 1, we consider the inverse problem of recovering
the double functions $p$ and $q$ from the given three spectra.

\section{main results}

Let $X_k(x,\la)$ and $Z_k(x,\la),$ $k=1,2,$ be the solutions of Eq.(\ref{1})
with the initial conditions
$$
X_1(0,\la)=X_2'(0,\la)=Z_1(T,\la)=Z_2'(T,\la)=1,
$$
$$
X_1'(0,\la)=X_2(0,\la)=Z_1'(T,\la)=Z_2(T,\la)=0.
$$
Denote by $L_0$ the boundary value problem (BVP) for Eq.(\ref{1}) with the conditions
$$
U_1(y)=U_2(y)=0,
$$
and $\om(\la):=\det[U_j(X_k)]_{j,k=1,2}$, and assume that $\om(\la)\not\equiv 0.$
The function $\om(\la)$ is an entire function of exponential type with order $1,$ and its zeros
$\Xi:=\{\xi_n\}_{n\in \mathbb{Z}}$ (counting multiplicities) coincide with the eigenvalues of $L_0$. The function $\om(\la)$
is called the characteristic function for $L_0$.

Denote $V_j(y):=y^{(j-1)}(T),\; j=1,2.$ Consider the BVP $L_j$, $j=1,2,$ for
Eq.(1) with the conditions $U_j(y)=V_1(y)=0.$ The eigenvalue sets $\Lambda_j:=
\{\la_{nj}\}_{n\in \mathbb{Z}}$ (counting multiplicities) of the BVP $L_j$ coincide with the zeros of the
characteristic function $\Delta_j(\la):=\det[U_j(X_k), V_1(X_k)]_{k=1,2}$.

For $\lambda\neq \lambda_{n1}$, let $\Phi(x,\la)$ be the solution of Eq.(\ref{1}) under the conditions $U_1(\Phi)=1$ and
$V_1(\Phi)=0.$ Denote Weyl-type function $M(\la):=U_2(\Phi).$ It is known \cite{BuYu} that for differential pencils with
classical two-point separated boundary conditions, the specification of the
Weyl function uniquely determines the double functions $p(x)$ and $q(x)$. In particular, in \cite{BuYu}
it is proved that differential pencils with
classical two-point separated boundary conditions is uniquely determined by specifying its Weyl
function, which is equivalent to specification of the spectra of two boundary value problems
with one common boundary condition, and a constructive procedure for
solving the inverse problem is given. However, in the case with
nonlocal conditions, it is not true; the specification of the Weyl-type
function $M(\la)$ does not uniquely determine the functions $p(x)$ and $q(x)$ (see counterexamples
in Section 5). For the nonlocal conditions the inverse problem is formulated as follows.

Throughout this paper the functions $\sigma_j(t)$ and the value $\int_0^T p(x)dx$ are known a priori.
And condition $S$: $\Lambda_1\cap\Xi=\emptyset$.

\smallskip
{\bf Inverse problem 1.} Given $M(\la)$ and $\om(\la),$ construct the functions $p(x)$ and $q(x).$

Let us formulate a uniqueness theorem for Inverse problem 1. For this purpose,
together with $(p, q)$ we consider another $(\tilde p,\tilde q),$ and we agree that
if a certain symbol $\al$ denotes an object related to $(p, q)$, then $\tilde\al$
will denote an analogous object related to $(\tilde p,\tilde q)$.

{\bf Theorem 1. }{\it Let condition $S$ be true. If $M(\la)=\tilde M(\la)$
and $\om(\la)=\tilde\om(\la),$ then $p(x)=\tilde p(x)$ and $q(x)\stackrel{a.e.}{=}
\tilde q(x)$ on $(0,T).$}

Thus, under condition $S,$ the specification $M(\la)$ and $\om(\la)$ uniquely
determines the function pairs $(p,q)$. We
note that if condition $S$ does not hold, then the specification $M(\la)$ and
$\om(\la)$ does not uniquely determine the functions (see Example 1 in
Section 5).

Consider the BVP $L_{11}$ for Eq.(\ref{1}) with the conditions $U_1(y)=V_2(y)=0.$
The eigenvalue set $\Lambda_{11}:=\{\la_{n1}^1\}_{n\in \mathbb{Z}}$ of the BVP $L_{11}$
coincide with the zeros of the characteristic function $\Delta_{11}(\la):=
\det[U_1(X_k), V_2(X_k)]_{k=1,2}$. Clearly, $\{\lambda_{n1}\}_{n\in \mathbb{Z}}\cap\{\lambda_{n1}^1\}_{n\in \mathbb{Z}}=\emptyset.$

{\bf Inverse problem 2. } Given $\{\la_{n1}, \la_{n1}^1\}_{n\in \mathbb{Z}}$,
construct $p(x)$ and $q(x).$

This inverse problem is a generalization of the well-known Borg's inverse
problem \cite{16} for Sturm-Liouville operators with classical two-point
separated boundary conditions.

{\bf Theorem 2. }{\it If $\la_{n1}=\tilde\la_{n1}, \la_{n1}^1=\tilde\la_{n1}^1$,
$n\in \mathbb{Z},$ then $p(x)=\tilde p(x)$ and $q(x)\stackrel{a.e.}{=}\tilde q(x)$ on $(0,T).$}

\section{Auxiliary Lemmas}

Denote $\Lambda^{\pm}:=\{\lambda:\pm\mbox{Im}\lambda\geq 0\}$. It is known (see, for example, \cite{4,Yurko1}) that there exists a fundamental system of
solutions $\{Y_k(x,\lambda)\}_{k=1,2}$ of Eq.(1) such that for $\lambda\in\Lambda^{\pm}$, $|\lambda|\to\iy$, $\nu=0,1$:
\begin{equation}\label{4}
\begin{array}{l}
Y_1^{(\nu)}(x,\lambda)\!=\!(i\lambda)^{\nu}\exp\left(i\left(\lambda x-\int_0^xp(t)dt\right)\right)(1\!+\!O(\lambda^{-1})),
\\
Y_2^{(\nu)}(x,\lambda)\!=\!(-i\lambda)^{\nu}\exp\left(-i\left(\lambda x-\int_0^xp(t)dt\right)\right)(1\!+\!O(\lambda^{-1})),
\end{array}
\end{equation}
and
\begin{align}\label{5}
\det[Y_k^{(\nu-1)}(x,\lambda)]_{k,\nu=1,2}=-2i\lambda(1+O(\lambda^{-1})).
\end{align}

{\bf Lemma 1. }(See \cite{yy}) {\it Let $\{W_k(x,\la)\}_{k=1,2}$ be a fundamental system of
solutions of Eq.(\ref{1}), and let $Q_j(y),\; j=1,2,$ be linear forms. Then}
\begin{align}\label{6}
\det[Q_j(W_k)]_{k,j=1,2}=\det[Q_j(X_k)]_{k,j=1,2}
\det[W_k^{(\nu-1)}(x,\la)]_{k,\nu=1,2}\,.
\end{align}

It follows from (\ref{5})-(\ref{6}) that
\begin{align}\label{7}
\det[Q_j(Z_k)]_{k,j=1,2}=\det[Q_j(X_k)]_{k,j=1,2},
\end{align}
and
\begin{align}\label{8}
\det[Q_j(Y_k)]_{k,j=1,2}=-2i\lambda(1+O(\lambda^{-1})) \det[Q_j(X_k)]_{k,j=1,2}\,.
\end{align}

Introduce the functions
$$
\vv(x,\la)=-\det[X_k(x,\la), U_1(X_k)]_{k=1,2},\;
\theta(x,\la)=\det[X_k(x,\la), U_2(X_k)]_{k=1,2},
$$
$$
\psi(x,\la)=\det[X_k(x,\la), V_1(X_k)]_{k=1,2}.
$$
Then
$$
U_1(\vv)=0,\;U_2(\vv)=\om(\la),\;V_1(\vv)=\Delta_1(\la),\;
V_2(\vv)=\Delta_{11}(\la),
$$
$$
U_1(\theta)=\om(\la),\; U_2(\theta)=0,\; V_1(\theta)=-\Delta_2(\la),
$$
$$
U_j(\psi)=\Delta_j(\la),\; V_1(\psi)=0,\; V_2(\psi)=-1.
$$
Moreover, by calculation, Eqs.(\ref{6})-(\ref{7}) yield that
\begin{equation}\label{9}
\begin{array}{l}
\det[\theta^{(\nu-1)}(x,\la), \vv^{(\nu-1)}(x,\la)]_{\nu=1,2}\!=\!\om(\la),\\
\det[\psi^{(\nu-1)}(x,\la), \vv^{(\nu-1)}(x,\la)]_{\nu=1,2}\!=\!\Delta_1(\la),
\end{array}
\end{equation}
\begin{align}\label{10}
\Delta_1(\la)=-U_1(Z_2),\;\Delta_2(\la)=-U_2(Z_2),\;\Delta_{11}(\la)=U_1(Z_1).
\end{align}
Note that the functions $\Phi, \psi, \vv$ and $\theta$ are all
the solutions of Eq.(\ref{1}) with some conditions, comparing boundary conditions on $\Phi, \psi, \vv$ and $\theta,$ we arrive at
\begin{align}\label{11}
\Phi(x,\la)=\frac{\psi(x,\la)}{\Delta_1(\la)},
\end{align}
\begin{align}\label{12}
\Phi(x,\la)=\frac{1}{\om(\la)}\,\Big(\theta(x,\la)+\frac{\Delta_2(\la)}{\Delta_1(\la)}\vv(x,\la)\Big).
\end{align}
Hence,
\begin{align}\label{13}
M(\la):=U_2(\Phi)=\frac{\Delta_2(\la)}{\Delta_1(\la)},
\end{align}
\begin{align}\label{14}
\det[\Phi^{(\nu-1)}(x,\la), \vv^{(\nu-1)}(x,\la)]_{\nu=1,2}=1.
\end{align}

Let $v_1(x,\la)$ and $v_2(x,\la)$ be the solutions of Eq.(\ref{1}) with
the conditions
$$
v_1(T,\la)=v_2'(T,\la)=1,\quad v_1'(T,\la)=0,\quad U_1(v_2)=0.
$$
Obviously,
\begin{equation}\label{15}
\begin{array}{l}
v_1(x,\la)=Z_1(x,\la),\; v_2(x,\la)=Z_2(x,\la)+N(\la)Z_1(x,\la),\\
\det[v_k^{(\nu-1)}(x,\la)]_{k,\nu=1,2}=1,
\end{array}
\end{equation}
where
\begin{align}\label{16}
N(\la)=\frac{\Delta_1(\la)}{\Delta_{11}(\la)}=-\frac{U_1(Z_2)}{U_1(Z_1)}\,.
\end{align}

Denote
$$
U_1^a(y):=\int_0^a y(t)d\sigma_1(t),\quad a\in(0,T].
$$
Clearly, $U_1=U_1^T,$ and if $\sigma_1(t)\equiv C$ (constant) for $t\ge a,$
then $U_1=U_1^a.$

For sufficiently small $\de>0,$ we denote
$$
\Pi_\de:=\{\lambda:\; \mbox{arg}\,\lambda\in[\de,\pi-\de]\},\ \
G_\de:=\{\lambda:\; |\lambda-\la_{n1}|\ge\de,\;\;\forall n\in \mathbb{Z}\},
$$
and
$$
G_\de':=\{\lambda:\; |\lambda-\la_{n1}^1|\ge\de,\;\;\forall n\in \mathbb{Z}\},
$$
where $\la_{n1}\in \Lambda_1$ and $\la_{n1}^1\in \Lambda_{11}$.

{\bf Lemma 2. }{\it For $\lambda\in\Pi_\de,\; |\lambda|\to\iy,$ we have
\begin{align}\label{17}
\Phi^{(\nu)}(x,\la)=
\frac{(i\lambda)^\nu}{H_1}\,\exp\left(i\left(\lambda x-\int_0^xp(t)dt\right)\right)(1+o(1)),\;x\in[0,T),
\end{align}
\begin{align}\label{18}
v_1^{(\nu)}(x,\la)=\frac{(i\lambda)^\nu}{2}
\,\exp\left(-i\left(\lambda(T-x)-\int_0^{T-x}p(t)dt\right)\right)(1+O(\lambda^{-1})),\; x\in[0,T),
\end{align}
\begin{equation}\label{19}
\begin{array}{l}
\Delta_1(\la)=-\frac{H_1}{2i\lambda}\,\exp\left(-i\left(\lambda T-\int_0^Tp(t)dt\right)\right)(1+o(1)),\\
\Delta_{11}(\la)=\frac{H_1}{2}\,\exp\left(-i\left(\lambda T-\int_0^T p(t)dt\right)\right)(1+o(1)).
\end{array}
\end{equation}

Let $\sigma_1(t)\equiv C$ (constant) for $t\ge a$ (i.e. $U_1=U_1^a$). Then for
$\lambda\in\Pi_\de,\; |\lambda|\to\iy,$}
\begin{align}\label{20}
&\vv^{(\nu)}(x,\la)=\frac{H_1}{2}\,(-i\lambda)^{\nu-1}
\exp\left(-i\left(\lambda x-\int_0^x p(t)dt\right)\right)\nonumber
\\
&\qquad\qquad\times[1+o(1)+O(\exp(i\lambda(2x-a)))],\; x\in(0,T],
\end{align}
\begin{align}\label{21}
&v_2^{(\nu)}(x,\la)=(-i\lambda)^{\nu-1}
\exp\left(i\left(\lambda(T-x)-\int_0^{T-x}p(t)dt\right)\right)\nonumber
\\
&\qquad\qquad\times[1+o(1)+O(\exp(i\lambda(2x-a)))],\; x\in[0,T).
\end{align}

{\it Proof.} The function $\Phi(x,\la)$ can be expressed as
\begin{align}\label{22}
\Phi(x,\la)=A_1(\la)Y_1(x,\lambda)+A_2(\la)Y_2(x,\lambda),
\end{align}
together with $U_1(\Phi)=1$ and $V_1(\Phi)=0,$ which yields that
\begin{align}\label{23}
A_1(\la)U_1(Y_1)+A_2(\la)U_1(Y_2)=1,\; A_1(\la)V_1(Y_1)+A_2(\la)V_1(Y_2)=0.
\end{align}
Using (4), one gets that for $\lambda\in\Pi_\de,\; |\lambda|\to\iy$:
\begin{align}\label{24}
U_1(Y_1)=H_1(1+o(1)),\; U_1(Y_2)=O(\exp(-i\lambda T)),
\end{align}
\begin{equation}\label{25}
\begin{array}{l}
V_1(Y_1)=\exp\left(i\left(\lambda T-\int_0^T p(t)dt\right)\right)(1+O(\lambda^{-1})),\\
V_1(Y_2)=\exp\left(-i\left(\lambda T-\int_0^T p(t)dt\right)\right)(1+O(\lambda^{-1})).
\end{array}
\end{equation}
Solving linear algebraic system (\ref{23}) by using (\ref{24})-(\ref{25}), we obtain
$$
A_1(\lambda)=H_1^{-1}(1+o(1)),\; A_2(\lambda)=O\left(\exp\left(2i\left(\lambda T-\int_0^T p(t)dt\right)\right)\right).
$$
Substituting these relations into (\ref{22}), we have proved (\ref{17}).
Formulas (\ref{18})-(\ref{21}) can be proved similarly, and are omitted. \B

By the well-known method (see, for example, \cite{4}) the
following estimates hold for $x\in(0,T),\; \lambda\in \Lambda^+:$
\begin{align}\label{26}
v_1^{(\nu)}(x,\la)=O\left(\lambda^{\nu}\exp\left(-i\left(\lambda(T-x)-\int_0^{T-x}p(t)dt\right)\right)\right),
\end{align}
\begin{align}\label{27}
\Phi^{(\nu)}(x,\la)=O\left(\lambda^{\nu}\exp\left(i\left(\lambda x-\int_0^xp(t)dt\right)\right)\right),\quad \rho\in G_\de.
\end{align}
Moreover, if $\sigma_1(t)\equiv C$ (constant) for $t\ge a$ (i.e. $U_1=U_1^a$),
then for $x\ge a/2,\; \lambda\in \Lambda^+$:
\begin{align}\label{28}
\vv^{(\nu)}(x,\la)=O\left(\lambda^{\nu-1}\exp\left(-i\left(\lambda x-\int_0^xp(t)dt\right)\right)\right),
\end{align}
\begin{align}\label{29}
v_2^{(\nu)}(x,\la)=O\left(\lambda^{\nu-1}\exp\left(i\left(\lambda(T-x)-\int_0^{T-x}p(t)dt\right)\right)\right),\quad \rho\in G'_\de.
\end{align}

\section{Proofs of Theorems}

{\textbf{Proof of Theorem 2}}
We know that the characteristic function $\Delta_1(\la)$ of the BVP $L_1$ is an entire function of order one with respect to $\lambda$.
Following the theory of Hadamard's factorization (see \cite{AHL}), $\Delta_1(\lambda)$ can be expressed as an
infinite product as
$$
\Delta_1(\lambda)=c_1e^{a_1\lambda}\prod_{n\in \mathbb{Z}}\left(1-\frac{\lambda}{\lambda_{n1}}\right)
e^{\frac{\lambda}{\lambda_{n1}}+\frac{1}{2}(\frac{\lambda}{\lambda_{n1}})^2+\cdots+\frac{1}{p}(\frac{\lambda}{\lambda_{n1}})^p},
$$
where $\lambda_{n1}$ are the eigenvalues of the problem $L_1$, $p$ is the genus of $\Delta_1(\lambda)$,
$c_1$ and $a_1$ are constants. Since for
the order $\rho$ of $\Delta_1(\lambda)$, $p\leq \rho\leq p+1$ (see \cite{AHL}), and $\Delta_1(\lambda)$ is
an entire function of exponential type with order $1,$ we find that the genus of $\Delta_1(\lambda)$ is $0$ or $1$ (that is, $p=0\vee 1$). Thus
$\Delta_1(\lambda)$ can be rewritten by
$$
\Delta_1(\lambda)=c_1e^{a_1\lambda}\prod_{n\in \mathbb{Z}}\left(1-\frac{\lambda}{\lambda_{n1}}\right)
e^{\frac{\lambda}{\lambda_{n1}}p}.
$$

Since $\Delta_1(\lambda)$ and $\tilde{\Delta}_1(\lambda)$ are both entire functions of order one with respect
to $\lambda$, and $\lambda_{n1}=\tilde{\lambda}_{n1}$ for all $n\in \mathbb{Z}$, by the Hadamard's factorization theorem, we may suppose
(the case when $\Delta_1(0)=0$ requires minor modifications)
$$
\Delta_1(\lambda)=c_1e^{a_1 \lambda}\prod_{n\in \mathbb{Z}}\left(1-\frac{\lambda}{\lambda_{n1}}\right)e^{\frac{\lambda}{\lambda_{n1}}p}
$$
and
$$
\tilde{\Delta}_1(\lambda)=\tilde{c}_1e^{\tilde{a}_1 \lambda}\prod_{n\in \mathbb{Z}}\left(1-\frac{\lambda}{\lambda_{n1}}\right)e^{\frac{\lambda}{\lambda_{n1}}\tilde{p}},
$$
for some constants $c_1,\tilde{c}_1,a_1,\tilde{a}_1$ and $p,\tilde{p}$, which can be determined from
the asymptotics.
From this we get for all $\lambda\in \mathbb{C}$
$$
\frac{\tilde{\Delta}_1(\lambda)}{\Delta_1(\lambda)}=\frac{\tilde{c}_1}{c_1}
e^{\left[(\tilde{a}_1-a_1)+(\tilde{p}-p)\sum_{n\in \mathbb{Z}}\frac{1}{\lambda_{n1}}\right]\lambda}.
$$
The expression (\ref{19}) implies that
$$
\Delta_1(\la)=-\frac{H_1}{2i\lambda}\,\exp\left(-i\left(\lambda T-\int_0^Tp(t)dt\right)\right)(1+o(1))
$$
and
$$
\tilde{\Delta}_1(\la)=-\frac{H_1}{2i\lambda}\,\exp\left(-i\left(\lambda T-\int_0^T\tilde{p}(t)dt\right)\right)(1+o(1)).
$$
We get that, using the assumption that $\int_0^Tp(t)dt=\int_0^T\tilde{p}(t)dt$,
$$
\frac{\tilde{\Delta}_1(\la}{\Delta_1(\la)}=1+o(1)\equiv\frac{\tilde{c}_1}{c_1}
e^{\left[(\tilde{a}_1-a_1)+(\tilde{p}-p)\sum_{n\in \mathbb{Z}}\frac{1}{\lambda_{n1}}\right]\lambda},
$$
which yields that
$$
(\tilde{a}_1-a_1)+(\tilde{p}-p)\sum_{n\in \mathbb{Z}}\frac{1}{\lambda_{n1}}=0,\ \tilde{c}_1=c_1.
$$
Consequently,
$\Delta_1(\la)\equiv\tilde\Delta_1(\la).$ Analogously, from $\lambda_{n1}^1=\tilde{\lambda}_{n1}^1$ for all $n\in \mathbb{Z}$ we get
$\Delta_{11}(\la)\equiv\tilde\Delta_{11}(\la).$ By virtue of (\ref{16}),
this yields
\begin{align}\label{30}
N(\la)\equiv\tilde N(\la).
\end{align}

Define the functions
\begin{equation}\label{31}
\begin{array}{l}
P_1(x,\la):=v_1(x,\la)\tilde v'_2(x,\la)\!-\!\tilde v'_1(x,\la)v_2(x,\la),\\
P_2(x,\la):=v_2(x,\la)\tilde v_1(x,\la)\!-\!\tilde v_2(x,\la)v_1(x,\la).
\end{array}
\end{equation}
Using (\ref{15}) and (\ref{30}), one gets
\begin{align*}
P_1(x,\la)&=(Z_1(x,\la)\tilde Z'_2(x,\la)-\tilde Z'_1(x,\la)Z_2(x,\la))
+(\tilde N(\la)-N(\la))Z_1(x,\la)\tilde Z'_1(x,\la)\\
&=Z_1(x,\la)\tilde Z'_2(x,\la)-\tilde Z'_1(x,\la)Z_2(x,\la),
\end{align*}
\begin{align*}
P_2(x,\la)&=Z_2(x,\la)\tilde Z_1(x,\la)-\tilde Z_2(x,\la)Z_1(x,\la)
+(N(\la)-\tilde N(\la))Z_1(x,\la)\tilde Z_1(x,\la)\\
&=Z_2(x,\la)\tilde Z_1(x,\la)-\tilde Z_2(x,\la)Z_1(x,\la).
\end{align*}
Thus, for each fixed $x\in (0,T),$ the functions $P_k(x,\la),\; k=1,2,$
are entire in $\la.$ On the other hand, taking (\ref{18}) and (\ref{21}) into
account we calculate for each fixed $x\ge T/2$ and $k=1,2:$
$$
P_k(x,\la)-\de_{1k}\Omega_1(x)=o(1),\; |\rho|\to\iy,\;\rho\in\Pi_\de,
$$
where $\de_{1k}$ is the Kronecker symbol and
$$
\Omega_1(x)=\frac{\exp[i\int_0^{T-x}(p(t)-\tilde{p}(t))dt]+\exp[-i\int_0^{T-x}(p(t)-\tilde{p}(t))dt]}{2}.
$$
Also, applying
(\ref{26}) and (\ref{29}), we get for $k=1,2$,
$$
P_k(x,\la)=O(1),\; |\rho|\to\iy,\;\rho\in G_\de'.
$$
Using the maximum modulus principle and Liouville's theorem for
entire functions, we conclude that
$$
P_1(x,\la)\equiv \Omega_1(x),\quad P_2(x,\la)\equiv 0,\quad x\ge T/2.
$$
From (\ref{31}) it yields
$$
v_1(x,\lambda)=\Omega_1(x)\tilde{v}_1(x,\lambda),\ v_2(x,\lambda)=\Omega_1(x)\tilde{v}_2(x,\lambda).
$$
Again, using the asymptotic expression (\ref{18}) for $v_1(x,\lambda)$ and $\tilde{v}_1(x,\lambda)$, we have
$$
\exp\left(-i(\lambda(T-x)-\int_0^{T-x}p(t)dt)\right)[1+o(1)]
$$
$$
\qquad\qquad=\Omega_1(x)\exp\left(-i(\lambda(T-x)-\int_0^{T-x}\tilde{p}(t)dt)\right)[1+o(1)],
$$
which leads to
$$
\exp\left(i\int_0^{T-x}(p(t)-\tilde{p}(t))dt\right)=\Omega_1(x).
$$
This deduces for $x\geq T/2$,
$$
\int_0^{T-x}(p(t)-\tilde{p}(t))dt=0,
$$
which yields
$$
p(x)=\tilde{p}(x) \mbox{ for } x\in \left[0,\frac{T}{2}\right].
$$
At this case we have $\int_0^{T-x}(p(t)-\tilde{p}(t))dt=0$ for $x\ge T/2$. Thus $\Omega_1(x)=1$ for $x\ge T/2$.

Together with (\ref{31}) this yields that for $x\ge T/2$,
\begin{align}\label{32}
v_k(x,\la)=\tilde v_k(x,\la),\;
Z_k(x,\la)=\tilde Z_k(x,\la),\; p(x)=\tilde p(x),\; q(x)\stackrel{a.e.}{=}\tilde q(x).
\end{align}

Next let us now consider the BVPs $L_1^a$ and $L_{11}^a$ for Eq.(\ref{1}) on
the interval $(0,T)$ with the conditions $U_1^a(y)=V_1(y)=0$ and
$U_1^a(y)=V_2(y)=0,$ respectively. Then, according to Eq.(\ref{10}), the
functions $\Delta_1^a(\la):=-U_1^a(Z_2)$ and $\Delta_{11}^a(\la)
:=U_1^a(Z_1)$ are the characteristic functions of $L_1^a$ and
$L_{11}^a,$ respectively. And
$$
U_1^{a/2}(Z_k)=U_1^{a}(Z_k)-\int_{a/2}^a Z_k(t,\la)d\sigma_1(t),
\quad k=1,2,
$$
hence
\begin{equation}\label{33}
\begin{array}{l}
\Delta_{1}^{a/2}(\la)=
\Delta_{1}^a(\la)+\int_{a/2}^a Z_2(t,\la)d\sigma_1(t),\\
\Delta_{11}^{a/2}(\la)=
\Delta_{11}^a(\la)-\int_{a/2}^a Z_1(t,\la)d\sigma_1(t).
\end{array}
\end{equation}
Let us use (\ref{33}) for $a=T.$ Since $\Delta_{1}^T(\la)=\Delta_{1}(\la),\;
\Delta_{11}^T(\la)=\Delta_{11}(\la),$ it follows from (\ref{32})-(\ref{33}) that
$$
\Delta_{1}^{T/2}(\la)=\tilde\Delta_{1}^{T/2}(\la),\quad
\Delta_{11}^{T/2}(\la)=\tilde\Delta_{11}^{T/2}(\la).
$$
Repeating preceding arguments subsequently for $a=T/2, T/4, T/8,\ldots,$
we conclude that $p(x)=\tilde{p}(x)$ and $q(x)\stackrel{a.e.}{=}\tilde q(x)$ on $(0,T).$ Theorem 2 is proved.\B

{\textbf{Proof of Theorem 1}}

Define the functions
\begin{equation}\label{34}
\begin{array}{l}
R_1(x,\la):=\Phi(x,\la)\tilde\vv'(x,\la)\!-\!\tilde\Phi'(x,\la)\vv(x,\la),\\
R_2(x,\la):=\vv(x,\la)\tilde\Phi(x,\la)\!-\!\tilde\vv(x,\la)\Phi(x,\la).
\end{array}
\end{equation}

Since $\Lambda_1\cap\Xi=\emptyset$ we can infer that $\Lambda_1\cap\Lambda_2=\emptyset$. Otherwise,
if a certain $\lambda\in \Lambda_1\cap\Lambda_2$ then $\lambda\in \Xi$. Thus $\lambda\in \Lambda_1\cap\Xi$;
this leads to a contradiction to the assumption that $\Lambda_1\cap\Xi=\emptyset$. Moreover, Eqs.
$M(\la)=\tilde M(\la)$, $M(\la)=\frac{\Delta_2(\lambda)}{\Delta_1(\lambda)}, \mbox{ and} \ \tilde{M}(\la)=\frac{\tilde{\Delta}_2(\lambda)}{\tilde{\Delta}_1(\lambda)}$
imply that
$$
\Delta_1(\lambda)=\tilde{\Delta}_1(\lambda),\ \Delta_2(\lambda)=\tilde{\Delta}_2(\lambda).
$$
It follows from (\ref{11}) and (\ref{34}) that
$$
R_1(x,\la)=\frac{1}{\Delta_1(\la)}
\Big(\psi(x,\la)\tilde\vv'(x,\la)-\tilde\psi'(x,\la)\vv(x,\la)\Big),
$$
$$
R_2(x,\la)=\frac{1}{\Delta_1(\la)}
\Big(\vv(x,\la)\tilde\psi(x,\la)-\tilde\vv(x,\la)\psi(x,\la)\Big).
$$
The above equations imply that for each fixed $x,$ the functions $R_k(x,\la)$ are
meromorphic in $\la$ with possible poles only at $\la=\la_{n1}$.
On the other hand, taking (\ref{12}) into account, we also get
\begin{align}\label{35}
R_1(x,\la)=\frac{1}{\om(\la)}
\Big(\theta(x,\la)\tilde\vv'(x,\la)-\tilde\theta'(x,\la)\vv(x,\la)\Big),
\end{align}
\begin{align}\label{36}
R_2(x,\la)=\frac{1}{\om(\la)}
\Big(\vv(x,\la)\tilde\theta(x,\la)-\tilde\vv(x,\la)\theta(x,\la)\Big).
\end{align}
The assumption that $\Lambda_1\cap\Xi=\emptyset$ tells us that the functions $R_k(x,\la)$ are regular at $\la=\la_{n1}$. Thus,
for each fixed $x,$ the functions $R_k(x,\la)$ are entire in $\la.$
Using (\ref{17}) and (\ref{20}), we can obtain for $x\ge T/2:$
$$
R_k(x,\la)-\de_{1k}\Omega_2(x)=o(1),\quad |\rho|\to\iy,\; \rho\in\Pi_\de,
$$
where
$$
\Omega_2(x)=\frac{\exp[i\int_0^x(p(t)-\tilde{p}(t))dt]+\exp[-i\int_0^x(p(t)-\tilde{p}(t))dt]}{2}.
$$
Also, using (27)-(28), we obtain for $x\geq T/2:$
$$
R_k(x,\la)=O(1),\quad |\rho|\to\iy,\; \rho\in G_\de.
$$
Therefore, $R_1(x,\la)\equiv \Omega_2(x),\; R_2(x,\la)\equiv 0$ for $x\geq T/2$.
Using the asymptotic expression (\ref{20}) for $\varphi(x,\lambda)$ and $\tilde{\varphi}(x,\lambda)$, we have
$$
\exp\left(-i(\lambda x-\int_0^xp(t)dt)\right)[1+o(1)]
$$
$$
\qquad\qquad=\Omega_2(x)\exp\left(-i(\lambda x-\int_0^x\tilde{p}(t)dt)\right)[1+o(1)].
$$
Thus
$$
\exp\left(i\int_0^x(p(t)-\tilde{p}(t))dt\right)=\Omega_2(x),
$$
which deduces for $x\geq T/2$,
$$
\int_0^x(p(t)-\tilde{p}(t))dt=0.
$$
this yields
$$
p(x)=\tilde{p}(x) \mbox{ for } x\in \left[\frac{T}{2},T\right].
$$
At this case we have $\int_0^x(p(t)-\tilde{p}(t))dt=0$ for $x\ge T/2$. Thus $\Omega_2(x)=1$ for $x\ge T/2$.
Together with (\ref{14}) and (\ref{34}), it yields
$$
\vv(x,\la)=\tilde\vv(x,\la),\; \psi(x,\la)=\tilde\psi(x,\la),\; p(x)=\tilde{p}(x),\;
q(x)\stackrel{a.e.}{=}\tilde q(x),\; x\ge T/2.
$$
Also, we obtain
$$
Z_k(x,\la)=\tilde Z_k(x,\la),\quad k=1,2,\quad x\ge T/2.
$$
Since
$$
\vv(x,\la)=U_1(Z_1)Z_2(x,\la)-U_1(Z_2)Z_1(x,\la)
$$
and
$$
\tilde{\vv}(x,\la)=U_1(\tilde{Z}_1)\tilde{Z}_2(x,\la)-U_1(\tilde{Z}_2)\tilde{Z}_1(x,\la)
$$
we have
$$
\vv(x,\la)=\Delta_{11}(\lambda)Z_2(x,\la)+\Delta_{1}(\lambda)Z_1(x,\la)
$$
and
$$
\tilde{\vv}(x,\la)=\tilde{\Delta}_{11}(\lambda)\tilde{Z}_2(x,\la)+\tilde{\Delta}_{1}(\lambda)\tilde{Z}_1(x,\la).
$$
Taking $x=T$ we get
$$
\vv(T,\la)=\Delta_{1}(\lambda),\ \tilde{\vv}(T,\la)=\tilde{\Delta}_{1}(\lambda),\ \vv'(T,\la)=\Delta_{11}(\lambda),\ \tilde{\vv}'(T,\la)=\tilde{\Delta}_{11}(\lambda).
$$
It follows from $\varphi(x,\la)=\tilde\varphi(x,\la)$ for $x\geq T/2$ that
$$
\Delta_{1}(\la)=\tilde\Delta_{1}(\la),\;\Delta_{11}(\la)=\tilde\Delta_{11}(\la).
$$
Using Theorem 2, we conclude that $p(x)=\tilde{p}(x)$ and $q(x)\stackrel{a.e.}{=}\tilde q(x)$ on $(0,T).$
Theorem 1 is proved.
\B

\section{Counterexamples}

{\textbf{Example 1} \ (To illustrate that if condition $S$ does not hold then Theorem 1 is false)}

Suppose that $T=\pi,\;U_1(y)=y(0),\;U_2(y)=y(\pi/2),$ $p(x)=p(x+\pi/2)$ and $q(x)=q(x+\pi/2)$ for $x\in(0,\pi/2),$
$p(x)\not\equiv p(\pi-x)$ and $q(x)\not\equiv q(\pi-x)$ for $x\in (0,\pi)$.

Take $\tilde p(x):=p(\pi-x)$ and $\tilde q(x):=q(\pi-x)$ for $x\in (0,\pi).$ Then BVP
$\tilde{L}_1$:
Eq.(\ref{1}) with  $\tilde{p}(x)=p(\pi-x)$ and $\tilde{q}(x)=q(\pi-x)$, $\tilde{y}(x)=y(\pi-x)$, and the conditions $U_1(\tilde{y})=V_1(\tilde{y})=0$;

BVP
$\tilde{L}_2$:
Eq.(\ref{1}) with $\tilde{p}(x)=p(\pi-x)$ and $\tilde{q}(x)=q(\pi-x)$, $\tilde{y}(x)=y(\frac{3\pi}{2}-x)$, and the conditions $U_2(\tilde{y})=V_1(\tilde{y})=0$;

BVP
$\tilde{L}_0$:
Eq.(\ref{1}) with $\tilde{p}(x)=p(\pi/2-x)$ (also equal to $p(\pi-x)$) and $\tilde{q}(x)=q(\pi/2-x)$ (also equal to $q(\pi-x)$), $\tilde{y}(x)=y(\frac{\pi}{2}-x)$, and the conditions $U_1(\tilde{y})=U_2(\tilde{y})=0$.

Here $\Delta_1(\lambda)$ is the characteristic function for Eq.(\ref{1}) with
$y(0)=0=y(\pi)$; $\Delta_2(\lambda)$ is the characteristic function for Eq.(\ref{1}) with
$y(\pi/2)=0=y(\pi)$;
$\omega(\lambda)$ is the characteristic function for Eq.(\ref{1}) with $y(0)=0=y(\pi/2)$.
From the above fact the following relations are true:
$$
\Delta_1(\la)=\tilde\Delta_1(\la),\; \Delta_2(\la)=\tilde\Delta_2(\la),\;
\om(\la)=\tilde\om(\la),
$$
and, in view of (\ref{13}), $M(\la)=\tilde M(\la).$

Note that $\Lambda_1$, $\Lambda_2$, and $\Xi$ are sets of zeros for characteristic functions
$\Delta_1(\lambda)$, $\Delta_2(\lambda)$ and $\omega(\lambda)$, respectively. Since $p(x)=p(x+\pi/2)$ and $q(x)=q(x+\pi/2)$ for $x\in(0,\pi/2),$
there holds
$\omega(\lambda)=\Delta_2(\lambda)$, i.e. $\Xi=\Lambda_2$. Thus for all $\lambda\in \Xi(=\Lambda_2)$ then it yields
$\lambda\in \Lambda_1$, which implies that $\Xi\cap\Lambda_1\neq \emptyset$.

Now $M(\la)=\tilde M(\la)$ and $\omega(\lambda)=\tilde{\omega}(\lambda)$, but $\Xi\cap\Lambda_1\neq \emptyset$. In Theorem 1
condition $S$ does not hold. In fact, at this case $p(x)\neq \tilde{p}(x):=p(\pi-x)$ and $q(x)\neq \tilde{q}(x):=q(\pi-x)$.
This means, that the specification of $M(\la)$ and $\om(\la)$ does not
uniquely determine the functions $p(x)$ and $q(x).$

{\textbf{Example 2} \ (To illustrate that even if condition $S$ and $M(\lambda)=\tilde{M}(\lambda)$ hold without
the assumption that $\om(\la)=\tilde\om(\la)$ then Theorem 1 is false)}

Suppose that $T=\pi,\;U_1(y)=y(0),\;U_2(y)=y(\pi-\al),$ where $\al\in(0,\pi/2).$

Let $p(x)\not\equiv p(\pi-x)$ and $q(x)\not\equiv q(\pi-x),$ and $(p(x),q(x))\equiv (0,0)$ for $x\in[0,\al_0]\cup
[\pi-\al_0,\pi],$ where $\al_0\in(0,\pi/2).$ If $\al<\al_0$, then $\la_{n2}=
\pi n/\al,\; n\in \mathbb{Z}\setminus\{0\}.$ Choose a sufficiently small $\al<\al_0$ such that
$\Lambda_1\cap\Lambda_2=\emptyset.$ Clearly, such choice is possible. Then
$\Lambda_1\cap\Xi=\emptyset,$ i.e. condition $S$ holds. Otherwise, if a certain $\lambda^*\in \Lambda_1\cap\Xi$, then
$\lambda^*\in \Lambda_1\cap\Lambda_2$; this contradicts to the fact that $\Lambda_1\cap\Lambda_2=\emptyset.$

Take $\tilde p(x):=p(\pi-x)$ and
$\tilde q(x):=q(\pi-x).$
Note that $\Delta_2(\lambda)$ is the characteristic function for the problem
$$
-y''(x)=\lambda^2 y(x), \ y(\pi-\alpha)=0=y(\pi);
$$
and
$\tilde{\Delta}_2(\lambda)$ is the characteristic function for the problem
$$
-y''(x)=\lambda^2 y(x), \ y(0)=0=y(\alpha).
$$
A simple calculation shows $\Lambda_2=\tilde{\Lambda}_2=\left\{\frac{\pi n}{\alpha}, n\in \mathbb{Z}\setminus\{0\}\right\}$.

At this case $\Delta_1(\la)=\tilde\Delta_1(\la),$
$\Delta_2(\la)=\tilde\Delta_2(\la),$ and consequently, $M(\la)=\tilde M(\la).$
Now in Theorem 1 condition $S$ holds, and $M(\la)=\tilde M(\la).$ In fact, in this example, $p(x)\neq \tilde{p}(x):=p(\pi-x)$ and $q(x)\neq \tilde{q}(x):=q(\pi-x)$.
So the true condition $S$ and the specification of $M(\la)$ does not uniquely
determine the functions $p(x)$ and $q(x).$

\section{Inverse problem from three spectra}
Fix $a\in(0,T).$ Consider Inverse problem 1 in the case when
$U_1(y):=y(0),\; U_2(y):=y(a).$ Then the boundary value problems $L_0, L_1, L_2$
take the forms
$$
L'_0:\ \mbox{Eq.}(\ref{1}) \mbox{ with } y(0)=y(a)=0,
$$
$$
L'_1:\ \mbox{Eq.}(\ref{1}) \mbox{ with } y(0)=y(T)=0,
$$
$$
L'_2:\ \mbox{Eq.}(\ref{1}) \mbox{ with } y(a)=y(T)=0.
$$

Denote by $\Lambda'_j=\{\la'_{nj}\}$ the spectrum of $L'_j\ (j=0,1,2),$ and assume that
$\Lambda'_0\cap\Lambda'_1=\emptyset$ (condition $S'$).

\smallskip
{\bf Inverse problem 4. } Given three spectra $\Lambda'_0, \Lambda'_1$ and
$\Lambda'_2$, construct $p(x)$ and $q(x).$

\smallskip
The following theorem is a consequence of Theorem 1.

\smallskip
{\bf Theorem 4. }{\it Let condition $S'$ hold. If $\Lambda'_j=\tilde\Lambda'_j$,
$j=0,1,2,$ then $p(x)=\tilde p(x)$ and $q(x)\stackrel{a.e.}{=}\tilde q(x)$ on $(0,T).$}

\smallskip
Comparing Theorem 4 with Theorem 1 we note that $\Lambda_0', \Lambda_1'$ and $\Lambda_2'$ correspond to
$\Xi,\Lambda_1$, and $\Lambda_2$ in Theorem 1, respectively.
In particular, Inverse problem 4 with $p(x)\equiv 0$ on $[0,T]$ was studied by many authors
(see, for example, \cite{17,18}).

\noindent {\bf Acknowledgments.}
The research work of the second author was supported by the Russian Ministry of
Education and Science (Grant 1.1436.2014K), and by Grant
13-01-00134 of Russian Foundation for Basic Research. The first author was supported in
part by the National Natural Science Foundation of China (11171152) and
Natural Science Foundation of Jiangsu Province of China (BK 20141392).


\begin{thebibliography}{99}

\bibitem{d} Day W.A. Extensions of a property of the heat equation to linear thermoelasticity
and order theories. Quart. Appl. Math. {\bf 40} (1982), 319-330.

\bibitem{g} Gordeziani N. On some nonlocal problems of the theory of elasticity. Bulletin
of TICMI {\bf 4} (2000), 43-46.

\bibitem{i} Ionkin N.I. The solution of a certain boundary value problem of the theory of
heat conduction with a nonclassical boundary condition. Differ. Equ. {\bf 13} (1997),
294-304. (in Russian)

\bibitem{1} Bitsadze A.V. and Samarskii A.A. Some elementary generalizations of linear
elliptic boundary value problems. Dokl. Akad. Nauk SSSR {\bf 185} (1969), 739-740.

\bibitem{yin} Yin Y.F. On nonlinear parabolic equations with nonlocal boundary conditions.
Journal of Mathematical Analysis and Applications {\bf 185} (1994), 161-174.

\bibitem{n} Nakhushev A.M. Equations of Mathematical Biology. Vysshaya Shkola, Moscow,
1995. (in Russian)

\bibitem{s} Schuegerl K. Bioreaction Engineering. Reactions Involving Microorganisms and
Cells, volume 1. John Wiley and Sons, 1987.

\bibitem{2} Marchenko V.A. Sturm-Liouville Operators and their Applications.
     Naukova Dumka,  Kiev, 1977;  English  transl., Birkh\"auser, 1986.
\bibitem{3} Levitan B.M. Inverse Sturm-Liouville Problems. Nauka,
     Moscow, 1984; English transl., VNU Sci. Press, Utrecht, 1987.
\bibitem{4} Freiling G. and Yurko V.A. Inverse Sturm-Liouville
     Problems and their Applications. NOVA Science Publishers, New York, 2001.
\bibitem{5} Yurko V.A. Method of Spectral Mappings in the Inverse Problem Theory.
     Inverse and Ill-posed Problems Series. VSP, Utrecht, 2002.
\bibitem{Yurko2} Yurko V. A. An inverse problem for pencils of differential operators. Matem. Sbornik {\bf 191} (2000), 137-160;
English transl., Sbornik: Mathematics {\bf 191} (2000), 1561-1586.
\bibitem{GGU} Gasymov M. G. and Gusejnov G. Sh. Determination of a diffusion operator from spectral data. Akad. Nauk
Azerb. SSR Dokl. {\bf 37} (1981), 19-23.
\bibitem{6} Yurko V.A. An inverse problem for integral operators. Matem. Zametki {\bf 37}
 (1985), 690-701 (Russian); English transl. in  Mathematical Notes {\bf 37} (1985), 378-385.
\bibitem{7} Yurko V.A. An inverse problem for integro-differential operators.
     Matem. Zametki {\bf 50} (1991), 134-146 (Russian); English transl. in
		 Mathematical Notes, {\bf 50} (1991), 1188-1197.
\bibitem{8} Kravchenko K.V. On differential operators with nonlocal boundary
     conditions. Differ. Uravn. {\bf 36} (2000), 464-469; English transl. in
		 Differ. Equations {\bf 36} (2000), 517-523.
\bibitem{9} Buterin S.A. The inverse problem of recovering the Volterra convolution
     operator from the incomplete spectrum of its rank-one perturbation.
		 Inverse Problems {\bf 22} (2006), 2223-2236.
\bibitem{10} Buterin S.A. On an inverse spectral problem for a convolution
     integro-differential operator. Results in Mathematics {\bf 50} (2007), 173-181.
\bibitem{11} Hryniv R., Nizhnik L.P., and Albeverio S. Inverse spectral problems for
     nonlocal Sturm-Liouville operators. Inverse Problems {\bf 23} (2007), 523-535.
\bibitem{12} Nizhnik L.P. Inverse nonlocal Sturm-Liouville problem. Inverse Problems
     {\bf 26} (2010), 125006 (9pp).
\bibitem{13} Kuryshova Y.V. and Chung-Tsun Shieh. Inverse nodal problem for
     integro-differential operators. Journal of Inverse and Ill-Posed Problems
		 {\bf 18} (2010), 357-369.
\bibitem{14} Freiling G. and Yurko V.A. Inverse problems for differential operators
     with a constant delay. Applied Mathematical Letters {\bf 25} (2012), 1999-2004.
\bibitem{15} Yang C. F. Trace and inverse problem of a discontinuous Sturm-Liouville
     operator with retarded argument. J. Math. Anal. Appl. {\bf 395} (2012), 30-41.
\bibitem{yy} Yurko V.A. and Yang C. F. Inverse problems for differential operators with nonlocal boundary conditions. Preprint, 2014.
\bibitem{BuYu} Buterin S. A. and Yurko V. A. Inverse spectral problem for pencils of differential operators on a finite interval.
Vestnik Bashkir. Univ. {\bf 4} (2006), 8-12.
\bibitem{16} Borg G. Eine Umkehrung der Sturm-Liouvilleschen Eigenwertaufgabe.
     Acta Math. {\bf 78} (1946), 1-96.
\bibitem{Yurko1} Yurko V. A. Recovering differential pencils on compact graphs. J. Differential Equations {\bf 244} (2008), 431-443.
\bibitem{AHL} Ahlfors L.~V. Complex Analysis, New York: McGraw-Hill, 1979.
\bibitem{17} Gesztesy F. and Simon B. On the determination of a potential from three
     spectra. Amer. Math. Soc. Transl. Ser.2, 189, 85-92, Amer. Math. Soc.,
		 Providence, RI, 1999.
\bibitem{18} Pivovarchik V. An inverse Sturm-Liouville problem by three spectra.
     Integral Equations Operator Theory {\bf 34} (1999), 234-243.
\end{thebibliography}
\end{document}